\documentclass{amsart}[12pt, leqno]
\usepackage{graphicx}

\theoremstyle{definition}

\usepackage{amscd,amssymb}

\begin{document}

\title[Varieties of associative algebras]
{Varieties of associative algebras\\
with an identity of third degree}
\author[L.A. Vladimirova, V.S. Drensky]{Lyubov A. Vladimirova, Vesselin S. Drensky}
\address{\hfill Received Dec. 7, 1983\newline
\phantom{X}
\newline Centre for Mathematics and Mechanics
\newline 1090 Sofia, P.O. Box 373\newline
\phantom{X}\newline
Current address: \newline Institute of Mathematics and Informatics\newline
Bulgarian Academy of Sciences\newline
Acad. G. Bonchev Str., Block 8, 1113 Sofia, Bulgaria}
\email{drensky@math.bas.bg}

\thanks{Published as Lyubov A. Vladimirova, Vesselin S. Drensky,
Varieties of associative algebras with an identity of third degree (Russian), Pliska Stud. Math. Bulgar. 8 (1986), 144-157.
\newline Zbl 0664.16015, MR0866654.}
\date{}


\begin{abstract}
We give a complete description of the varieties of associative algebras over a field of characteristic zero which satisfy a polynomial identity of third degree.
\end{abstract}

\maketitle

One of the main problems in the theory of PI-algebras is the description of all varieties of algebras with a prescribed property.
When one describes all subvarieties of a given variety, usually this is made in the language of lattices (see the survey by Artamonov \cite{2}).
In the case of associative algebras over a field of characteristic 0 the simplest non-trivial case is the case of varieties with an identity of third degree.
This problem was studied by Nagata \cite{12}, Klein \cite{10}, Olsson and Regev \cite{13}, Anan'in and Kemer \cite{1}, Regev \cite{14, 15, 16}, James \cite{9}.

The purpose of this paper is to complete these investigations and to give a complete description of the varieties with an identity of degree 3.
Clearly, this problem is equivalent to the description of all T-ideals containing such an identity.
In the proof we use the method which allowed one of the authors to obtain the lattice in the case of some varieties of unitary algebras and Lie algebras \cite{4, 5}.

{\bf  1. Preliminaries.}
We shall use notations similar to those in \cite{4, 5}: $K$ is a fixed field of characteristic 0;
$S_n$ is the symmetric group acting on the set $\{1,2,\ldots,n\}$;
$A_m=K\langle x_1,\ldots,x_m\rangle$ is the free non-unitary associative algebra with generators $x_1,x_2,\ldots, x_m$ (we shall use also other symbols for the free generators);
$A_m^{(n)}$ is the homogeneous component of degree $n$ of $A_m$;
$P_n$ is the set of multilinear polynomials in $A_m^{(n)}$,
$P_n$ is a left $S_n$-module with respect to the action $\sigma:x_{i_1}\cdots x_{i_n}\to x_{\sigma(i_1)}\cdots x_{\sigma(i_n)}$, $\sigma\in S_n$.
The irreducible $S_n$-module corresponding to the partition $\lambda=(\lambda_1,\ldots,\lambda_r)$ of $n$ will be denoted by $M(\lambda)=M(\lambda_1,\ldots,\lambda_r)$.
On the other hand, the algebra $A_m$ is isomorphic to the tensor algebra of the $m$-dimensional vector space
and has a natural structure of a left $GL(m,K)$-module. The corresponding irreducible modules will be $N_m(\lambda)=N_m(\lambda_1,\ldots,\lambda_r)$
(and $N_m(\lambda)=0$ if $\lambda_{m+1}>0$).

The lattice of the $GL(m,K)$-submodules of $A_m^{(n)}$ can be embedded into the lattice of the $S_n$-submodules of $P_n$. The embedding is realized in the following way.
There are canonical generators $f(x_1,\ldots,x_r)$ for $N_m(\lambda_1,\ldots,\lambda_r)\subset A_m^{(n)}$, $r\leq m$, and $e(x_1,\ldots,x_n)$ for $M(\lambda_1,\ldots,\lambda_r)\subset P_n$,
$\sum \lambda_i=n$. The linearization of $f(x_1,\ldots,x_r)$ gives $e(x_1,\ldots,x_n)$; The ``symmetrization'' of $e(x_1,\ldots,x_n)$ gives, up to a nonzero multiplicative constant, $f(x_1,\ldots,x_r)$
if $r\leq m$ and 0 if $r>m$.

Information on representation theory of $S_n$ and $GL(m,K)$ can be found e.g. in \cite{3} or \cite{7}, and in the language of polynomial identities in \cite{1} and \cite{4}.
The algebra $A_m$ has a natural multigrading: $A_m=\sum A_m^{(k_1,\ldots,k_m)}$, $\sum k_i\geq 1$,
where $A_m^{(k_1,\ldots,k_m)}$ is the vector space of homogeneous polynomials of degree $k_i$ in $x_i$.
This grading agrees with the action of $GL(m,K)$: if $g=\sum \xi_ie_{ii}$ is a diagonal linear transformation from $GL(m,K)$, then $A_m^{(k_1,\ldots,k_m)}$ is the set of eigenvectors
corresponding to the eigenvalue $\xi_1^{k_1}\cdots\xi_m^{k_m}$ of $g$.
Let $\mathfrak M$ be a variety of algebras with a T-ideal of its identities $T({\mathfrak M})$. The relatively free algebra $F_m({\mathfrak M})=A_m/(T({\mathfrak M})\cap A_m)$ and
the set of multilinear elements $P_n({\mathfrak M})=P_n/(T({\mathfrak M})\cap P_n)$ preserve the action of $GL(m,K)$ and $S_n$, respectively, and $F_m({\mathfrak M})$ and $P_n({\mathfrak M})$
have the same module structure \cite{4, 6}: If $P_n({\mathfrak M})=\sum\varkappa_{\lambda}M(\lambda)$, then $F_m({\mathfrak M})=\sum\varkappa_{\lambda}N_m(\lambda)$.
(All the sums of modules are direct and the subscripts in the notation for the vector subspaces of $F_m({\mathfrak M})$ are the same as in the case of $A_m$.)
For our purposes it is more convenient to work with the representations of $GL(m,K)$, but following the traditions, the results will be stated in the language of representations of $S_n$.
The Hilbert series of the graded vector space $F_m({\mathfrak M})$
\[
H(F_m({\mathfrak M}),t_1,\ldots,t_m)=\sum\dim F_m^{(k_1,\ldots,k_m)}({\mathfrak M})t_1^{k_1}\cdots t_m^{k_m}
\]
is equal to $\sum\varkappa_{\lambda}S_{\lambda}(t_1,\ldots,t_m)$, where $S_{\lambda}(t_1,\ldots,t_m)$ is the Schur function. As it was noted by Berele \cite{6}, this series completely determines
the module structure of $F_m({\mathfrak M})$.

{\bf Lemma 1.1.} {\it Let $S_{\lambda}(t_1,\ldots,t_m)=\sum\beta_{\lambda}^{(k)}t_1^{k_1}\cdots t_m^{k_m}$, where the summation is on all $k=(k_1,\ldots,k_m)$.
Let us order lexicographically all partitions of $n$. If $P_n({\mathfrak M})=\sum\varkappa_{\lambda}M(\lambda)$, then}
\begin{equation}\label{(1)}
\varkappa_{\lambda}=\dim A_r^{(\lambda)}-\dim(T({\mathfrak M})\cap A_r^{(\lambda)})-\sum_{\mu>\lambda}\varkappa_{\mu}\beta_{\mu}^{(\lambda)}.
\end{equation}

The proof follows from the fact that the matrix of the coefficients $\beta_{\mu}^{(\lambda)}$ is unitriangular \cite{11}.

For direct applications of Lemma 1.1 we shall need several coefficients of the Schur functions. They are given in a table in \cite{11}.

{\bf Lemma 1.2.} $\beta_{(4)}^{(3,1)}=\beta_{(3,1)}^{(2,2)}=\beta_{(4)}^{(2,1,1)}=\beta_{(2,2)}^{(2,1,1)}=1$, $\beta_{(3,1)}^{(2,1,1)}=2$,

\noindent $\beta_{(3,2)}^{(3,1,1)}=\beta_{(3,2)}^{(2,2,1)}=1$, $\beta_{(3,1,1)}^{(2,2,1)}=\beta_{(2,2,1)}^{(2,1,1,1)}=2$,
$\beta_{(3,2)}^{(2,1,1,1)}=\beta_{(3,1,1)}^{(2,1,1,1)}=3$.

For a variety $\mathfrak M$ with an identity of degree 3 at least one of the irreducible components of $P_3$ is contained in $T({\mathfrak M})$.
It is known that
\[
P_3\cong KS_3\cong M(3)+2M(2,1)+M(1^3).
\]
The generators of the corresponding submodules are the linearizations of the polynomials
$x^3$ for $M(3)$, $\alpha[x,y]x+\beta x[x,y]$, $\alpha\not=0$ or $\beta\not=0$ for $M(2,1)$ and the standard identity $S_3(x_1,x_2,x_3)$ for $M(1^3)$.
We shall consider separately these three cases. In each of the cases our considerations are in the following way:
1) We determine the module structure of $P_n({\mathfrak M})$. Here we use essentially earlier results from \cite{1, 9, 12, 14, 15, 16}.
2) For every irreducible component of $P_n({\mathfrak M})$ we find the consequences of higher degree.
It is sufficient to find only the consequences in $P_{n+1}({\mathfrak M})$.
We shall denote this graphically as follows: The irreducible modules will be denoted by dots and by arrows we shall denote their consequences.
This immediately gives the complete description of the T-ideals with identity of degree three and the lattice of the subvarieties of the corresponding varieties.

In the sequel all computations will be in the relatively free algebra of the variety $\mathfrak M$.

{\bf 2. The variety $\mathfrak M$ defined by the identity}
\begin{equation}\label{(2)}
\alpha[x,y]x+\beta x[x,y]=0.
\end{equation}
The structure of $P_n({\mathfrak M})$ was almost completely described in \cite{1}.
Depending on the coefficients $\alpha$ and $\beta$ the problem is reduced to four cases:

\noindent
2.1. $\alpha\not=0$, $\beta\not=0$, $\alpha-\beta\not=0$, $\alpha+\beta\not=0$ (the general case);

\noindent 2.2. $\alpha=0$ (the case $\beta=0$ is handled in a similar way);

\noindent 2.3. $\alpha-\beta=0$;

\noindent $\alpha+\beta=0$.

{\bf The case 2.1.} $\alpha\beta(\alpha-\beta)(\alpha+\beta)\not=0$.

{\bf Proposition 2.1.1.} {\it $P_n({\mathfrak M})\cong M(n)$, $n=1$ or $n\geq 4$, $P_2({\mathfrak M})\cong M(2)+M(1^2)$,
$P_3({\mathfrak M})\cong M(3)+M(2,1)+M(1^3)$.

The corresponding modules are generated by the linearizations of $x^n$ for $M(n)$, $[x_1,x_2]x_1$ for $M(2,1)$, $S_n(x_1,\ldots,x_n)$ for $M(1^n)$, $n=2,3$.}

{\bf Proof.} For $n\leq 3$ the modules $P_n({\mathfrak M})$ are the same in all four cases and can be described directly using that $P_n\cong KS_n$. Let $n\geq 4$.
It follows from \cite{1} that if the module $P_n({\mathfrak M})$ is nonzero, then it is generated by
$h_n(x_1,\ldots,x_n)=\sum x_{\sigma(1)}\cdots x_{\sigma(n)}$. The algebra of (commuting) polynomials $K[x]$ belongs to $\mathfrak M$ and $x^n\not=0$ in $K[x]$.
The polynomial $h_n(x_1,\ldots,x_n)$ is the linearization of $x^n$ and therefore $P_n({\mathfrak M})\cong M(n)$.

{\bf Theorem 2.1.2.} {\it Let the linearization of the identity $f(x_1,\ldots,x_r)$ generate an irreducible $S_n$-submodule of $P_n({\mathfrak M})$. Then
the consequences of higher degree of $f(x_1,\ldots,x_r)$ are equivalent to the identities:

{\rm (a)} $x_1x_2\cdots x_{n+1}$ if $f=x^n$;

{\rm(b)} $[x_1,x_2]x_3$ if $f=[x_1,x_2]$;

{\rm(c)} $[x_1,x_2]x_1$ and $S_3(x_1,x_2,x_3)$ do not have consequences of higher degree} (see Fig. 1).
\begin{center}
\includegraphics[width=4.6cm]{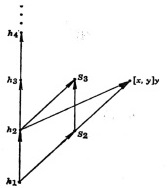}\\
Fig. 1
\end{center}

{\bf Proof.} For $n>3$ the polynomial $x^n$ is equivalent in $\mathfrak M$ to $x_1x_2\cdots x_n$. The statements (a) and (b) can be verified directly.
The case (c) follows from the fact that in $K[x]$ we have $[x_1,x_2]x_1=S_3(x_1,x_2,x_3)=0$ and $x^n\not=0$.

{\bf The case 2.2.} $\alpha=0$. Then the variety $\mathfrak M$ is defined by the identity $x[x,y]=0$.

{\bf Proposition 2.2.1.} {\it $P_n({\mathfrak M})\cong M(n)+M(n-1,1)$, $n\geq 4$.
The submodules $M(n)$ and $M(n-1,1)$ are generated by the linearizations of $x^n$ and $f_{n-1}=[x,y]x^{n-2}$, respectively.}

{\bf Proof.} It was shown in \cite{1} that for $n\geq 4$ $P_n({\mathfrak M})$ is a sum of not more than two irreducible submodules.
The linearization of $[x,y]x^{n-2}$ generates $M(n-1,1)$ in $P_n({\mathfrak M})$. Hence it is sufficient to show that $[x,y]x^{n-2}\not=0$ in $F({\mathfrak M})$.
Let
\[
A_1=\langle a,b\mid a^2=a,b^2=ab=0,ba=b\rangle
\]
be a two-dimensional associative algebra. It is isomorphic to the algebra of $2\times 2$ matrices
$\{k_1e_{11}+k_2e_{21}\mid k_1,k_2\in K\}$. Then $A_1\in\mathfrak M$ and $[a,b]a^{n-2}=-b\not=0$.

{\bf Theorem 2.2.2.} {\it Let $f(x_1,\ldots,x_r)\in F({\mathfrak M})$. Then the consequences of higher degree of $f$ are equivalent to the identities:

{\rm (a)} $x_1x_2\cdots x_{n+1}$ if $f=x^n$;

{\rm (b)} $[x_1,x_2]x_1^{n-1}$ if $f=[x_1,x_2]x_1^{n-2}$, $n>2$;

{\rm (c)} $S_3(x_1,x_2,x_3)$ does not have consequences of higher degree} (see Fig. 2).

\begin{center}
\includegraphics[width=4.6cm]{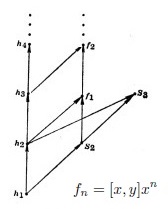}\\
Fig. 2
\end{center}

{\bf Proof.} (a) It is sufficient to show that in $F({\mathfrak M})$ the identities $x_{n+1}$ and $[x_1,x_2]x_1^{n-1}$ follow from $x^n$.
The identity $xyx=x^2y$ follows from the identity $x[x,y]=0$ which defines the variety $\mathfrak M$. But
$[x,y]x^{n-1}=xyx^{n-1}-yx^n=xyx^n-1$ and $[x,y]x^{n-1}=xyx^{n-1}=x^2yx^{n-2}+\cdots=x^ny$.

The case (b) is obvious.

(c) The algebra $A_1$ from the proof of Proposition 2.2.1 is two-dimensional and $S_3(x_1,x_2,x_3)=0$ in $A_1$. Additionally, $A_1\in\mathfrak M$, $x_1^4\not=0$, $[x_1,x_2]x_1^2\not=0$ in $A_1$.
Therefore $S_3(x_1,x_2,x_3)$ does not have consequences of higher degree.

{\bf The case 2.3.} $\alpha=\beta$. The variety $\mathfrak M$ is defined by the identity
\begin{equation}\label{(3)}
[x^2,y]=0
\end{equation}
and the linearization of (\ref{(3)}) gives
\begin{equation}\label{(4)}
f(x_1,x_2,y)=[x_1x_2+x_2x_1,y]=0
\end{equation}

{\bf Proposition 2.3.1.} {\it $P_4({\mathfrak M})\cong M(4)+M(2^2)$, $P_n({\mathfrak M})=M(n)$ for $n\geq 5$. The submodule $M(2^2)$ is generated by the linearization of $x[x,y]y$.}

{\bf Proof.} In was shown in \cite{1} that $P_4({\mathfrak M})$ is a sum of not more than two non-isomorphic irreducible submodules and $P_n({\mathfrak M})\cong M(n)$ for $n\geq 5$.
It follows from (\ref{(3)}) that $[x,yzt]=0$, $x^3y-xyx^2=x^2yx-yx^3=x^3y-yx^3=0$ and these three identities are linearly independent. Using (\ref{(1)}) and Lemma 1.2 we obtain
\[
0\leq \varkappa_{(3,1)}\leq \dim A_2^{(3,1)}-1\cdot \beta_4^{(3,1)}=0,\,\varkappa_{(3,1)}=0.
\]
It also follows from (\ref{(3)}) that
\begin{equation}\label{(5)}
x^2y^2-xy^2x=yx^2y-y^2x^2=x^2y^2-y^2x^2=xyxy-yxyx=0
\end{equation}
and the identities in (\ref{(5)}) are linearly independent. Hence
\begin{equation}\label{(6)}
\dim(T({\mathfrak M})\cap A_2^{(2,2)})\geq 4.
\end{equation}
All consequences of (\ref{(3)}) in $A_2^{(2,2)}$ are obtained from (\ref{(4)}) as linear combinations of
\begin{equation}\label{(7)}
f(u_1,u_2,u_3)u_4=f(u_1u_2,u_3,u_4)=f(u_1,u_2,u_3u_4)=0,
\end{equation}
where two of the variables $u_1,u_2,u_3,u_4$ are equal to $x$ and the other two -- to $y$.
A direct verification shows that all possible identities (\ref{(7)}) follow from (\ref{(5)}), i.e. we have an equality in (\ref{(6)}).

Lemmas 1.1 and 1.2 give that $\varkappa_{(2,2)}=1$. It is easy to check that the polynomial $[x,y]^2$
which in $A_2^{(4)}$ generates a submodule isomorphic to $N_2(2^2)$ is proportional to $x[x,y]y$ modulo the identity (\ref{(5)}).

{\bf Theorem 2.3.2} {\it Let $f(x_1,\ldots,x_r)\in F({\mathfrak M})$. Then the consequences of higher degree of $f$ are equivalent to the identities:

{\rm (a)} $x_1x_2\cdots x_{n+1}$ if $f=x^n$, $n\geq 1$;

{\rm (b)} $x_1[x_1,x_2]x_2$ if $f=[x_1,x_2]x_1$ or $f=S_3(x_1,x_2,x_3)$;

{\rm (c)} $x_1[x_1,x_2]x_2$ does not have consequences of higher degree} (see Fig. 3).

\begin{center}
\includegraphics[width=4.6cm]{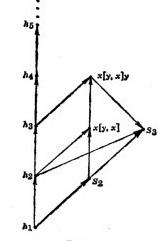}\\
Fig. 3
\end{center}

{\bf Proof.} We shall show that $x^3=0$ implies $x[x,y]y=0$, i.e. that $x_1x_2x_3x_4=0$ is a consequence of $x^3=0$. The partial linearization of $x^3=0$ gives
\[
x^2y+xyx+yx^2=0\text{ and hence }g(x,y)=x^2y^2+xy^2x+y^2x^2=0;
\]
\[
h_3(xy,y,x)-(g(x,y)+g(y,x))=yxyx-x^2y^2+2xyxy-2y^2x^2=0.
\]
We apply (\ref{(5)}) and obtain that $x[x,y]y=0$.

Similarly we obtain that $S_3(x,y,xy)=0$ is equivalent to the identity $x[x,y]y=0$.

The other cases of the theorem are trivial.

{\bf The case 2.4. $\alpha+\beta=0$.} Then the identity (\ref{(2)}) is equivalent to $[x,y,z]=0$.

{\bf Proposition 2.4.1.} (see \cite{1}.) {\it
\[
P_n({\mathfrak M})\cong M(n)+M(n-1,1)+M(n-2,1^2)+\cdots+M(1^n)
\]
and the modules $M(n-k,1^k)$ are generated by the linearizations of
\[
g_{k+1,n-k-1}=S_{k+1}(x_1,\ldots,x_{k+1})x_1^{n-k-1}.
\]
}

{\bf Theorem 2.4.2.} {\it All consequences of higher degree of $S_{k+1}(x_1,\ldots,x_{k+1})x_1^{n-k-1}$ in $F({\mathfrak M})$ are equivalent to
\[
\varepsilon S_k(x_1,\ldots,x_k)x_1^{n-k+1},S_{k+1}(x_1,\ldots,x_{k+1})x_1^{n-k},
\]
\[
S_{k+2}(x_1,\ldots,x_{k+2})x_1^{n-k-1},S_{k+3}(x_1,\ldots,x_{k+3})x_1^{n-k-2},
\]
where $\varepsilon=1$ if $k$ is even and $\varepsilon=0$ if $k$ is odd} (see Fig. 4).

\begin{center}
\includegraphics[width=8.6cm]{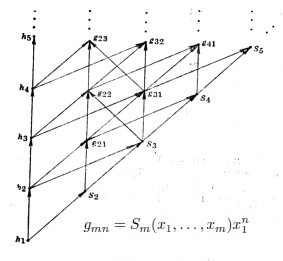}\\
Fig. 4
\end{center}

{\bf Proof.} We shall use the consequences of $[y,x,x]=0$
\begin{equation}\label{(8)}
[x,y]z=z[x,y];[x,y][z,t]=-[x,t][z,y].
\end{equation}
It follows from \cite[Proposition 1.1]{5} that all multilinear consequences of degree $n+1$ of the module $M(n-k,1^k)$ are among
\[
M(n-k+2,1^{k-1}),M(n-k+1,1^k),M(n-k,1^{k+1}),M(n-k-1,1^{k+2}).
\]
One can see directly that in the case of $M(n-k+2,1^{k-1})$, $k=2m$, the polynomial $S_{2m+1}(x_1,\ldots,x_{2m},x_1^2)x_1^{n-k-1}$ is proportional to
$S_{2m}(x_1,\ldots,x_{2m})x_1^{n-k+1}$.

When $k=2m-1$, then $S_{2m-1}(x_1,\ldots,x_{2m-1})x_1^{n-k+1}$ does not follow from
$S_{2m}(x_1,\ldots,x_{2m})x_1^{n-k+1}$ because the latter polynomial is a consequence of the product of $m$ commutators of length two
and the former polynomial cannot be written as a product of more than $m-1$ commutators.
 It follows from \cite{8} that
$S_{2k}(x_1,\ldots,x_{2k})$ and $S_{2k+1}(x_1,\ldots,x_{2k+1})$ are proportional, respectively, to
\[
[x_1,x_2][x_3,x_4]\cdots[x_{2k-1},x_{2k}]
\]
and
\[
[x_1,x_2][x_3,x_4]\cdots[x_{2k-1},x_{2k}]x_{2k+1}+[x_2,x_3][x_4,x_5]\cdots[x_{2k},x_{2k+1}]x_1
\]
\[
+\cdots +[x_{2k+1},x_1][x_2,x_3]\cdots[x_{2k-2},x_{2k-1}]x_{2k}.
\]

The case $M(n-k+1,1^k)$ is obvious. If $k$ is even, then the linearization of $S_{k+1}(x_1,\ldots,x_k,[x_{k+1},x_{k+2}])x_1^{n-k-1}$ is nonzero in $\mathfrak M$
and generates $M(n-k,1^{k+1})$. If $k$ is odd, then, taking into account \cite{8}, it is sufficient to multiply $S_{2m}x_1^{n-k-1}$ by $x_{2m+1},x_1,x_2,\ldots,x_{2m}$
and to take the sum of the products. The result will be nonzero and proportional to $S_{2m+1}(x_1,\ldots,x_{2m+1})x_1^{n-k-1}$.

In $S_{k+1}(x_1,\ldots,x_{k+1})x_1^{n-k-1}$ we replace $x_1$ by $x_1+[x_{k+2},x_{k+3}]$ and take the linear component in $x_{k+2}$. After some computations it will turn out
that it is equivalent to the identity $S_{k+3}(x_1,\ldots,x_{k+3})x_1^{n-k-2}$. This completes the proof of the theorem.

{\bf 3. The variety $\mathfrak M$ defined by the identity $S_3(x_1,x_2,x_3)=0$.}

{\bf Proposition 3.1.} {\it
\[
P_3({\mathfrak M})\cong M(3)+2M(2,1); P_4({\mathfrak M})\cong M(4)+2M(3,1)+M(2^2);
\]
\[
P_n({\mathfrak M})\cong M(n)+2M(n-1,1), n\geq 5.
\]
The irreducible components of $P_n({\mathfrak M})$ are generated by the linearizations of the polynomials:
\[
x^n\text{ for }M(n),f_n(x,y)=\alpha[x,y]x^{n-2}+\beta x^{n-2}[x,y]\text{ for }M(n-1,1), [x,y]^2\text{ for }M(2^2).
\]}

{\bf Proof.} The module structure of $P_n({\mathfrak M})$ was completely described in \cite{9}. For the proof of the proposition it is sufficient to show that
the polynomials $f_n(x,y)$ and $[x,y]^2$ are nonzero in $F_2({\mathfrak M})$. The module $N_2(2^2)$ participates with multiplicity 2 in the decomposition
of $A_2^{(4)}$ into a sum of irreducible $GL(2,K)$-modules. But
\[
S_3(x_1,x_2,[x_1,x_2])=\sum(-1)^{\sigma}[x_1,x_2,x_{\sigma(1)}]x_{\sigma(2)}+[x_1,x_2]^2
\]
and this is the only consequence of $S_3(x_1,x_2,x_3)$ which generates a module $N_2(2^2)$. Hence $[x_1,x_2]^2\not=0$ in $F_2({\mathfrak M})$.
The two-dimensional algebras (see the case 2.2)
\[
B_1=\{k_1e_{11}+k_2e_{21}\mid k_1,k_2\in K\}\text{ and }B_2=\{k_1e_{11}+k_2e_{12}\mid k_1,k_2\in K\}
\]
belong to $\mathfrak M$ and for them $f_n(x,y)\not=0$, $n\geq 3$. This completes the proof of the proposition.

{\bf Theorem 3.2.} {\it Let $f(x_1,\ldots,x_n)\in F({\mathfrak M})$. Then the consequences of higher degree of $f$ are equivalent to the identities:

{\rm (a)} $x_1x_2\cdots x_{n+1}$ if $f=x^n$;

{\rm (b)} $[x_1,x_2]x_3\cdots x_{n+1}$ if $f=[x,y]x^{n-2}$;

{\rm (c)} $x_1\cdots x_{n-1}[x_n,x_{n+1}]$ if $f=x^{n-2}[x,y]$;

{\rm (d)} $[x_1,x_2]x_3\cdots x_{n+1}$ and $x_1\cdots x_{n-1}[x_n,x_{n+1}]$ if
$f=\alpha[x,y]x^{n-2}+\beta x^{n-2}[x,y]$, $\alpha\not=0$, $\beta\not=0$;

{\rm (e)} $[x,y]^2$ does not have consequences of higher degree} (see Fig. 5).

\begin{center}
\includegraphics[width=8.6cm]{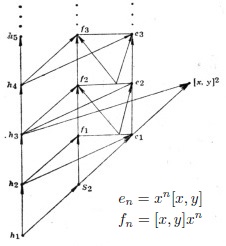}\\
Fig. 5
\end{center}

{\bf Proof.} (a) In $F({\mathfrak M})$
\[
S_3(x,y,x^2)=x[y,x]x=0.
\]
On the other hand, it follows from $x^n=0$
\[
0=x^{n-1}yx+x^{n-2}yx^2+\cdots+xyx^{n-1}=(n-1)x^{n-1}yx=(n-1)x^ny+x^{n-1}[y,x],
\]
i.e. $x^{n-1}[y,x]=0$ follows from $x^n=0$. In the same way we derive the consequence $[y,x]x^{n-1}=0$.
For $n\geq 4$ we obtain that all irreducible submodules in $P_{n+1}({\mathfrak M})$ are consequences of $x^n=0$.
Hence $x_1\cdots x_{n+1}=0$ follows from $x^n=0$ in $F({\mathfrak M})$.
For $n=3$ we have to show also that $[x,y]^2=0$. The identities
\[
x^2y^2+xy^2x+y^2x^2=y^2x^2+yx^2y+x^2y^2=0
\]
are consequences of $x^3=0$ and hence
\begin{equation}\label{(9)}
xy^2x=yx^2y.
\end{equation}
Additionally, $0=[x^2y+xyx+yx^2,y]$ and, by (\ref{(9)}),
\begin{equation}\label{(10)}
x^2y^2+xyxy=yxyx+y^2x^2.
\end{equation}
In $(x^2y+xyx+yx^2)x=x^2yx+xyx^2=0$ we replace $x$ by $x+y$ and take the homogeneous component of degree 2 in $x$:
\[
x^2y^2+2xy^2x+yxyx+y^2x^2+xyxy=0.
\]
In follows from (\ref{(9)}) and (\ref{(10)}) that
\[
x^2y^2+xy^2x+xyxy=0.
\]
Since $x^2y^2+xy^2x+y^2x^2=0$ we obtain $xyxy=y^2x^2$. Similarly $yxyx=x^2y^2$. Now
\[
S_3(x,y,xy)=[x,y]^2+x[y,x]y=0
\]
implies
\[
x[x,y]y=[x,y]^2=y[y,x]x,
\]
\[
2[x,y]^2=x[x,y]y+y[y,x]x=x^2y^2-yxyx+y^2x^2-xyxy=0.
\]
In this way we obtain that $[x,y]^2=0$ is a consequence of $x^3=0$ in $F({\mathfrak M})$.

(b) It is obvious that $[x,y]x^{n-1}=0$ follows from $[x,y]x^{n-2}=0$ in $F({\mathfrak M})$. For $n\geq 4$ this means that
$[x_1,x_2]x_3\cdots x_{n+1}=0$ On the other hand, $x^{n-1}[x,y]=0$ does not follow from $[x,y]x^{n-2}=0$ because
the latter identity holds for the algebra $B_2$ and the former does not hold.

The case $n=3$ was handled in the case 2.

(c) The proof is similar to that in (b).

(d) As a consequence of the identity
\[
f(x,y)=\alpha[x,y]x^{n-2}+\beta x^{n-2}[x,y]=0,\alpha\not=0,\beta\not=0,
\]
we obtain
\[
f(x,y)x=\alpha[x,y]x^{n-1}=0,xf(x,y)=\beta x^{n-1}[x,y]=0
\]
(since $x[x,y]x=0$ in $F({\mathfrak M})$). This gives the proof for $n\geq 4$. The case $n=3$ was considered in part 2.

(e) The identity $[x,y]^2=0$ does not have consequences in $F({\mathfrak M})$
because it is satisfied in the algebras $B_1$ and $B_2$ from the proof of Proposition 3.1
and the nonzero polynomials of degree 5 are not equal to 0 in these algebras.

{\bf 4. The variety $\mathfrak M$ defined by the identity $x^3=0$.}

The identity $x^3=0$ is equivalent to its linearization
\[
h_3(x_1,x_2,x_3)=\sum x_{\sigma(1)}x_{\sigma(2)}x_{\sigma(3)}=0,\sigma\in S_3.
\]
It is known \cite{12} that $x_1x_2x_3x_4x_5x_6=0$ in $F({\mathfrak M})$.
Hence for the description of the module structure of $P_n({\mathfrak M})$ it is sufficient to handle the cases $n=4,5$. First we shall prove:

{\bf Theorem 4.1.} $P_4({\mathfrak M})\cong M(3,1)+M(2^2)+2M(2,1^2)+M(1^4)$.

{\bf Proof.} It will be carried out in several steps and will follow from Lemmas 4.2 -- 4.5.

{\bf Lemma 4.2.} {\it The identity
\[
x^2yx+xyx^2=x^2y^2-yxyx=0
\]
holds in $F({\mathfrak M})$.}

{\bf Proof.} Since $\displaystyle 0=\frac{1}{2}h_3(x,x,y)x=x^2yx+xyx^2+yx^3$ in $F({\mathfrak M})$, we obtain $x^2yx+xyx^2=0$.
We also have
\[
\frac{1}{2}[h_3(x,x,y),y]=(x^2y^2+xyxy)-(yxyx+y^2x^2)=0.
\]
In $x^2yx+xyx^2=0$ we replace $x$ by $x+y$ and take the homogeneous component of degree 2 in $y$:
\[
x^2y^2+2xy^2x+yxyx+y^2x^2+xyxy=2(x^2y^2+xy^2x+xyxy)=0.
\]
Since $\displaystyle \frac{1}{2}h_3(x,x,y^2)=x^2y^2+xy^2x+y^2x^2=0$, we obtain $y^2x^2=xyxy$ and similarly $x^2y^2=yxyx$.

{\bf Lemma 4.3.} {\it The module $M(3,1)$ participates with multiplicity $1$ in the decomposition of $P_4({\mathfrak M})$
into a sum of irreducible submodules.}

{\bf Proof.} By (\ref{(1)}) and Lemma 1.2
\[
\varkappa_{(3,1)}=\dim A_2^{(3,1)}-\dim(T({\mathfrak M})\cap A_2^{(3,1)}).
\]
(The module $N(4)$ is generated by $x^4$ and hence is contained in $T({\mathfrak M})$.) Here $\dim A_2^{(3,1)}=4$ and
$\varkappa_{(3,1)}=4-\dim(T({\mathfrak M})\cap A_2^{(3,1)})$.

All homogeneous consequences of degree (3,1) of $x^3=0$ are obtained from $h_3(x,y,z)$ replacing $x,y,z$ by $x,x^2,y,xy,yx$
and by multiplication from the right by $x$ and $y$. One can check directly that all identities of degree (3,1) are linear combinations
of the identities $x^3y=yx^3=x^2yx-xyx^2=0$ which are linearly independent.
Hence $\dim(T({\mathfrak M})\cap A_2^{(3,1)})=3$ and $\varkappa_{(3,1)}=1$.

{\bf Lemma 4.4.} {\it The multiplicity of $M(2^2)$ in $P_4({\mathfrak M})$ is equal to $1$.}

{\bf Proof.} It follows from (\ref{(1)}) that
\[
\varkappa_{(2,2)}=\dim A_2^{(2,2)}-\dim(T({\mathfrak M})\cap A_2^{(2,2)})-1\cdot \beta_{(3,1)}^{(2,2)}=5-\dim(T({\mathfrak M})\cap A_2^{(2,2)}).
\]
As in the previous lemma one can see that the identities
\[
xyxy-y^2x^2=xy^2x+x^2y^2+y^2x^2=yxyx-x^2y^2=yx^2y+y^2x^2+x^2y^2=0
\]
are linearly independent and all homogeneous consequences of $x^3=0$ of degree (2,2) are their linear combinations. Hence
$\dim(T({\mathfrak M})\cap A_2^{(2,2)})=4$ and $\varkappa_{(2,2)}=1$.

{\bf Lemma 4.5.} {\it The multiplicity of $M(2,1^2)$ in $P_4({\mathfrak M})$ is equal to $2$.}

{\bf Proof.} It follows from (\ref{(1)}) that
\[
\varkappa_{(2,1^2)}=\dim A_3^{(2,1^2)}-\dim(T({\mathfrak M})\cap A_3^{(2,1^2)})-1\cdot\beta_{(3,1)}^{(2,1^2)}-1\cdot\beta_{(2,2)}^{(2,1^2)}.
\]
The vector space $A_3^{(2,1^2)}$ has a basis of 12 elements:
\[
xyxz,xyzx,yxzx,xzxy,xzyx,zxyx,yx^2z,yzx^2,
x^2yz,zx^2y,zyx^2,x^2zy.
\]
We take the three identities
\[
xyxz+x^2yz+yx^2z=xyzx+x^2yz+yzx^2=yxzx+yx^2z+yzx^2=0
\]
together with the similar identities obtained by changing the places of $y$ and $z$ and the identity
\[
x^2yz+yx^2z+yzx^2+x^2zy+zx^2y+zyx^2=0.
\]
The matrix with rows consisting of the coordinates of these identities with respect to the above basis has rank 7.
Hence these seven identities are linearly independent.
As in the previous lemmas, one can show that this is the maximal number of linearly independent polynomials in $T({\mathfrak M})\cap A_3^{(2,1^2)}$.
Therefore, by Lemma 1.2
\[
\varkappa_{(2,1^2)}=12-7-1\cdot 2-1\cdot 1=2.
\]
The module $M(1^4)$ participates in the decomposition of $P_4({\mathfrak M})$ because by \cite[Proposition 1.1]{5} $S_4(x,y,z,t)=0$
cannot be obtained as a consequence of $h_3(x,y,z)=0$.

The next theorem gives the module structure of $P_5({\mathfrak M})$. Its proof follows from Lemmas 4.7 -- 4.11.

{\bf Theorem 4.6.} $P_5({\mathfrak M})\cong M(3,2)+M(3,1^2)+M(2,1^3)$.

{\bf Lemma 4.7.} {\it The modules $M(5)$ and $M(4,1)$ do not participate in the decomposition of $P_5({\mathfrak M})$.}

{\bf Proof.} In virtue of Lemmas 1.1 and 1.2 it is sufficient to show that the monomials
$x^4y,x^3yx,x^2yx^2,xyx^3,yx^4$ are qual to zero in $F_2({\mathfrak M})$. This is not obvious only for $x^2yx^2$.
Since $x(x^2y+xyx+yx^2)x=0$, we obtain that $x^2yx^2=0$.

{\bf Lemma 4.8.} {\it The multiplicity of $M(3,2)$ in the decomposition of $P_5({\mathfrak M})$ into a sum of irreducible $S_5$-modules is equal to $1$.}

{\bf Proof.} By Lemmas 1.1 and 4.7
\[
\varkappa_{(3,2)}=\dim A_2^{(3,2)}-\dim(T({\mathfrak M})\cap A_2^{(3,2)}).
\]
Here $\dim A_2^{(3,2)}=10$. We multiply from the left or from the right by $x$ or by $y$ the identities
\[
x^2yx+xyx^2=x^2y^2+xy^2x+y^2x^2=x^2y^2+xy^2x+xyxy=y^2x^2-xyxy=x^2y^2-yxyx=0
\]
from the proof of Lemma 4.2 and obtain that
\[
xyxyx=0,x^2y^2x=-xy^2x^2=xyx^2y=-x^2yxy=yxyx^2=-yx^2yx.
\]
Together with the obvious identities $x^3y^2=y^2x^3=yx^3y=0$ we obtain 9 linearly independent identities of degree (3,2).
After a direct verification we can prove that all consequences of $x^3=0$ which are of degree 3 in $x$ and of degree 2 in $y$
are consequences of these 9 identities. Therefore $\dim(T({\mathfrak M})\cap A_2^{(3,2)})=4$ and $\varkappa_{(3,2)}=1$.

{\bf Lemma 4.9.} {\it There is one module $M(3,1^2)$ in the decomposition of $P_5({\mathfrak M})$.}

{\bf Proof.} By Lemmas 1.1, 1.2, 4.7 and 4.8
\[
\varkappa_{(3,1^2)}=\dim A_3^{(3,1^2)}-\dim(T({\mathfrak M})\cap A_3^{(3,1^2)})-1\cdot \beta_{(3,2)}^{(3,1^2)},
\]
where $\beta_{(3,2)}^{(3,1^2)}=1$ and $\dim A_3^{(3,1^2)}=20$.
From the identity in \cite{8}
\begin{equation}\label{(11)}
yx^2zt+tyx^2z=0
\end{equation}
for $t=x$ we obtain that $yx^2zx+xyx^2z=0$. From $xyx^2=-x^2yx$ we obtain $yxzx^2+x^2yxz=0$. Also
\[
h_3(x,x,yxz)=yxzx^2+xyxzx+x^2yxz=0
\]
which implies $xyxzx=0$. Similarly $xzxyx=0$. In $xyx^2+x^2yx=0$ we replace $y$ by $yz$ and obtain
$xyzx^2+x^2yzx=0$. We replace $x$ by $x+v$ in (\ref{(11)}) and take the multilinear component. Then for $y=x=v$ we obtain
\[
x^2vzx+xvxzx+x^3vz+x^2vxz=0,
\]
i.e. the identity $x^2yzx+x^2yxz=0$ holds in $\mathfrak M$. In this way we obtain 9 identities in $\mathfrak M$:
\[
x^3yz=yzx^3=yx^3z=xyxzx=0,
\]
\[
xyzx^2=-x^2yzx=x^2yxz=-xyx^2z=yx^2zx=-yxzx^2.
\]
We obtain 9 more identities exchanging the places of $y$ and $z$.
A direct verification shows that these 18 identities are linearly independent.
In order to show this, it is sufficient to write them with respect to any basis of $A_3^{(3,1^2)}$ and
to find the rank of the matrix of coordinates. The most convenient is to consider the following basis of $A_3^{(3,1^2)}$:
A factor $x^3$ participates in the first 6 monomials, then we continue with the monomials $xyxzx$ and $xzxyx$.
The next group of basis elements consists of $xyzx^2,x^2yzx, x^2yxz,xyx^2z,yx^2zx$ and of 5 more monomials obtained exchanging $y$ and $z$.
In the end we add the elements $yxzx^2$ and $zxyx^2$. Then the $18\times 18$ minor consisting of the first 18 columns of the matrix is equal to 1.
On the other hand, if we find all consequences of $x^3=0$ of degree $(3,1^2)$, it is easy to see that they are linear combinations
of these 18 identities. Therefore $\dim(T({\mathfrak M})\cap A_3^{(3,1^2)})=18$ and $\varkappa_{(3,1^2)}=1$.

{\bf Lemma 4.10.} {\it The module $M(2^2,1)$ does not participate in the decomposition of $P_5({\mathfrak M})$.}

{\bf Proof.} Lemmas 1.1, 1.2 and already established results in Section 4 give
\[
\varkappa_{(2^2,1)}=\dim A_3^{(2^2,1)}-\dim(T({\mathfrak M})\cap A_3^{(2^2,1)})-1\cdot\beta_{(3,2)}^{(2^2,1)}-1\cdot\beta_{(3,1^2)}^{(2^2,1)},
\]
\[
\dim A_3^{(2^2,1)}=30,\beta_{(3,2)}^{(2^2,1)}=2,\beta_{(3,1^2)}^{(2^2,1)}=1.
\]
It follows from the identity $x^2u+xux+ux^2=0$ that $xz(yxy+y^2x+xy^2)=0$, i.e.
\[
xzyxy=-xzxy^2-xzy^2x=(x^2z+zx^2)y^2+x^2(zy^2)+(zy^2)x^2.
\]
On the other hand,
\[
xzyxy=(-x^2(zy)-(zy)x^2)y=-x^2zy^2+z(y^2x^2+x^2y^2).
\]
Therefore
\[
2x^2zy^2+zx^2y^2+zy^2x^2=-x^2zy^2+zx^2y^2+zy^2x^2\text{ and }x^2zy^2=0.
\]
Similarly $y^2zx^2=0$. The obtained identities easily imply the following identities in $F({\mathfrak M})$:
$zxyxy=zy^2x^2$ (in Lemma 4.4 we showed that $xyxy=y^2x^2$)
\[
xyzxy=x^2y^2z+zx^2y^2,x^2yzy=-x^2y^2z,xy^2zx=-x^2y^2z,xyxyz=y^2x^2z,
\]
\[
xzy^2x=-zy^2x^2, xzyxy=zx^2y^2+zy^2x^2,xy^2xz=-x^2y^2z-y^2x^2z,
\]
\[
yxzxy=y^2x^2z+zx^2y^2, zxy^2x=-zx^2y^2-zy^2x^2,xyxzy=x^2y^2z+y^2x^2z,
\]
\[x^2y^2z+y^2x^2z+zx^2y^2+zy^2x^2=0, xzxy^2=-zx^2y^2.
\]
We obtain similar identities changing the places of $x$ and $y$. One can show that these 27 identities are linearly independent.
For this purpose it is sufficient to show that the matrix of the coefficients of these identities with respect to any basis of $A_3^{(2^2,1)}$
has rank 27. It is convenient to work with respect to the following basis of $A_3^{(2^2,1)}$: $x^2zy^2, y^2zx^2$,
followed by all monomials with the property that the two symbols $x$ (or two symbols $y$) are not neighbours,
$x^2y^2z,zx^2y^2,y^2x^2z,zy^2x^2$. This explains that $\dim(T({\mathfrak M})\cap A_3^{(2^2,1)})\geq 27$
and $0\leq \varkappa_{(2^2,1)}\leq 30-27-3=0$ and hence $\varkappa_{(2^2,1)}=0$.

{\bf Lemma 4.11.} {\it The module $M(2,1^3)$ participates with multiplicity $1$ in the decomposition of $P_5({\mathfrak M})$
into a sum of irreducible $S_5$-submodules.}

{\bf Proof.} By Lemmas 1.1, 1.2 and the already known multiplicities we have
\[
\varkappa_{(2,1^3)}=\dim A_4^{(2,1^3)}-\dim(T({\mathfrak M})\cap A_4^{(2,1^3)})-1\cdot \beta_{(3,2)}^{(2,1^3)}-1\cdot \beta_{(3,1^2)}^{(2,1^3)},
\]
\[
\dim A_4^{(2,1^3)}=60,\beta_{(3,2)}^{(2,1^3)}=\beta_{(3,1^2)}^{(2,1^3)}=3.
\]
Hence the problem is reduced to the determining the dimension of $T({\mathfrak M})\cap A_4^{(2,1^3)}$. It follows from the identity $x^2zy^2=0$
\begin{equation}\label{(12)}
x^2yzt+x^2ytz=yztx^2+zytx^2=0.
\end{equation}
It was established in \cite{8} that $\mathfrak M$ satisfies the identity
\begin{equation}\label{(13)}
yx^2zt+tyx^2z=0.
\end{equation}
In the linearization $xyztu+yxztu+xyzut+yxzut=0$ of $x^2zy^2=0$ we replace $u$ by $x$ and using $xvx=-x^2v-vx^2$, we obtain
\begin{equation}\label{(14)}
x^2yzt+yztx^2+yx^2zt+yzx^2t=0.
\end{equation}
It is easy to check also the identities
\begin{equation}\label{(15)}
\begin{split}
xyxzt=-x^2yzt-yx^2zt,\,yxzxt=-yx^2zt-yzx^2t,\\
xyzxt=-x^2yzt-yzx^2t,\,yxztx=-yx^2zt-yztx^2,\\
xyztx=-x^2yzt-yztx^2,\,yzxtx=-yzx^2t-yztx^2.
\end{split}
\end{equation}
Permuting the variables $y,z,t$ in the six identities (\ref{(15)}), we obtain 30 new identities. Similarly, we obtain 18 more identities from (\ref{(12)}) and (\ref{(13)}).
In this way we have 54 identities similar to (\ref{(12)}), (\ref{(13)}) and (\ref{(15)}). It has turned out that 53 of them are linearly independent
but all 54 are linearly dependent. In order to show this it is sufficient to express these polynomials with respect to the basis of $A_4^{(2,1^3)}$
consisting of the monomials
\[
xyxzt,xyzxt, xyzxt, yxzxt,yxztx, yzxtx,
\]
and 30 more monomials obtained by permuting of $y,z,t$. In the end we have 24 monomials where $x^2$ participates:
\[
\begin{split}
x^2yzt,x^2ytz,ytzx^2,tyzx^2,x^2tyz,x^2tzy,tzyx^2,ztyx^2,\\
x^2zty,x^2zyt,zytx^2,yx^2zt,yx^2tz,zx^2ty,zx^2yt,tx^2yz,\\
tx^2zy,yzx^2t,ytx^2z,ztx^2y,zyx^2t,tyx^2z,tzx^2y,yztx^2.
\end{split}
\]
The $60\times 54$ matrix of the coordinates of the considered identities with respect to this basis has rank 53.
As a result of listing all possible consequences of $x^3=0$ of degree $(2,1^3)$ we obtain that all consequences of the considered kind
are linear combinations of the identities like (\ref{(12)}), (\ref{(13)}), (\ref{(15)}). Hence
$\dim(T({\mathfrak M})\cap A_4^{(2,1^3)})=53$ and $\varkappa_{(2,1^3)}=1$.

For the proof of Theorem 4.6 it is sufficient to pay attention that
\[
S_5(x_1,x_2,x_3,x_4,x_5)=0
\]
is proportional to the identity
\[
\sum(-1)^{\sigma}h_3(x_{\sigma(1)}x_{\sigma(2)},x_{\sigma(3)}x_{\sigma(4)},x_{\sigma(5)})=0,
\]
and hence $S_5=0$ is a consequence of $x^3=0$.

{\bf Proposition 4.12.} {\it The linearizations of the following elements generate irreducible $S_n$-submodules in $P_n({\mathfrak M})$:
\[
M(1^3): S_3(x_1,x_2,x_3);\,M(2,1): \alpha x[x,y]+\beta[x,y]x, (\alpha,\beta)\not=(0,0);
\]
\[
M(1^4): S_4(x_1,x_2,x_3,x_4);\,M(2,1^2):\alpha f_1(x,y,z)+\beta f_2(x,y,z),
\]
\[
\text{where }f_1(x,y,z)=x^2[y,z]-[y,z]x^2,f_2(x,y,z)=yx^2z-zx^2y,(\alpha,\beta)\not=(0,0);
\]
\[
M(2^2):xy^2x;\, M(3,1): x^2yx;\,M(2,1^3):x^2S_3(y,z,t)+2d_3(y,z,t;x^2,1)
\]
\[
\text{(here }d_3(x_1,x_2,x_3;y_1,y_2)=\sum(-1)^{\sigma}x_{\sigma(1)}y_1x_{\sigma(2)}y_2x_{\sigma(3)}
\text{ is the Capelli identity)};
\]
\[
M(3,1^2):x[y,z]x^2;\,M(3,2):x^2y^2x.
\]
}

{\bf Proof.} The main idea of the proof of the proposition is the following. First, we fix a partition $\lambda=(\lambda_1,\ldots,\lambda_r)$ of $n$.
In the notation of \cite[\S 1]{4} we consider the Young tableau corresponding to the partition $\lambda$ for a suitable permutation $\tau$.
The symmetrization $f_{\tau}(x_1,\ldots,x_r)$ generates an irreducible $GL(m,K)$-submodule in $F_m({\mathfrak M})$ isomorphic to $N_m(\lambda)$.
Using already known consequences of $x^3=0$ we rework $f_{\tau}(x_1,\ldots,x_t)$ and show that in $F_m({\mathfrak M})$ it is equivalent to the identities in the statement of Proposition 4.12.

The cases $M(1^3)$, $M(2,1)$ and $M(1^4)$ are obvious. We shall assume that $x_1=x$, $x_2=y$, $x_3=z$, $x_4=t$.

For the partition $(2,1^2)$ and the identity permutation $\tau_1$ we obtain:
\[
f_{\tau_1}(x,y,z)=S_3(x,y,z)x=-x^2[y,z]+[y,z]x^2+yx^2z-zx^2y.
\]
Then for the permutation $\tau_2=(34)$ we have:
\[
f_{\tau_2}(x,y,z)=[x,y]xz+[y,z]x^2
+[z,x]xy=[y,z]x^2-x^2[y,z]+2(zx^2y-yx^2z).
\]
Using the identities in the proof of Lemma 4.5, we can show that the polynomials $f_{\tau_1}$ and $f_{\tau_2}$ are linearly independent in $F_3({\mathfrak M})$.
This means that their linearizations generate two different isomorphic submodules $M(2,1^2)$ in $P_4({\mathfrak M})$.
And it is easy to see that $f_{\tau_1}$ and $f_{\tau_2}$ are linear combinations of $f_1$ and $f_2$.

The case $(2^2)$. In this case we obtain for the identity permutation $\tau$
\[
f_{\tau}(x,y)=x^2y^2-yx^2y-xy^2x+y^2x^2.
\]
It follows from the proof of Lemma 4.4 that $f_{\tau}(x,y)=-3xy^2x\not=0$ in $F({\mathfrak M})$.

The cases $(3,1),(3,2),(3,1^2)$ can be handled in a similar way using the identity permutation and the identities from Lemmas 4.3, 4.8 and 4.9.

The case $(2,1^3)$. Using the identities (\ref{(12)}), (\ref{(13)}), (\ref{(14)}), (\ref{(15)}) it is easy to see that from all irreducible submodules corresponding to the standard Young tableaux,
only for the identity permutation we obtain a nonzero generator: $x^2S_3(y,z,t)+2d_3(y,z,t;x^2,1)$.

The next theorem gives the complete description of the T-ideals containing the polynomial $x^3$.

{\bf Theorem 4.13.} {\it In the variety $\mathfrak M$ defined by the identity $x^3=0$
all consequences of higher degree from the identities generating irreducible $S_n$-submodules in $P_n({\mathfrak M})$ are given in Fig. 6 (a), (b), (c).
In this figure
\[
f=yx^2z-zx^2y, e=x^2S_3(y,z,t)+2d_3(y,z,t;x^2,1),g=x^2[y,z]-[y,z]x^2.
\]}
\begin{center}
\includegraphics[width=8.6cm]{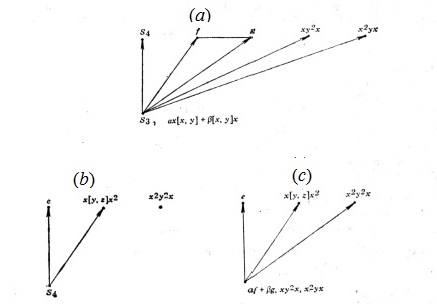}\\
Fig. 6
\end{center}

We shall complete the proof of the theorem in several steps:

1) In order to handle the case of Fig. 6 (a) it is sufficient to find the consequences of higher degree from $x^3$ in the varieties defined by the identities, respectively,
\[
S_3(x,y,z)=0,\alpha x[x,y]+\beta[x,y]x=0,(\alpha,\beta)\not=(0,0),
\]
(see Fig. 1, 2, 3, 4, 5).

2) In the proof of Proposition 4.12 we already mentioned that modulo the identities (\ref{(12)}), (\ref{(13)}), (\ref{(14)}), (\ref{(15)})
\[
S_4(x,y,z,t)x= x^2S_3(y,z,t)+2d_3(y,z,t;x^2,1),
\]
i.e. the element $x^2S_3+2d_3$ obviously is a consequence of $S_4(x,y,z,t)$ in $F({\mathfrak M})$.

In (\ref{(14)}) we take the alternating sum in $y,z,t$ and obtain
\begin{equation}\label{(16)}
x^2S_3(y,z,t)+S_3(y,z,t)x^2=0.
\end{equation}
Replacing $t$ by $x$ and using (\ref{(12)}) we see that the identity $x[y,z]x^2=0$ is a consequence of $S_4(x,y,z,t)$ in $F({\mathfrak M})$.
In this way, we obtain the case of Fig. 6 (b).
By \cite[Proposition 1.1]{5} the identity $x^2y^2x=0$ does not follow from $S_4(x,y,z,t)=0$.

3) We shall establish that
\[
x^2S_3(y,z,t)+2d_3(y,z,t;x^2,1)=0
\]
is a consequence of $xy^2x=0$ and $x^2yx=0$.

We take the alternating sum in (\ref{(11)}) and obtain
\begin{equation}\label{(17)}
d_3(y,z,t;x^2,1)+d_3(y,z,t;1,x^2)=0.
\end{equation}
In $yx^2y=0$ we replace $y$ by $y+z$ and then replace $z$ by $zx$. We compute the result using (\ref{(14)}) and (\ref{(15)}) and obtain
$ytx^2z+2ytzx^2=0$. This implies
\begin{equation}\label{(18)}
d_3(y,z,t;x^2,1)+2S_3(y,z,t)x^2=0.
\end{equation}
We linearize $x^2yzx$. Then it follows from (\ref{(15)}) that
\[
2tyx^2z=yzx^2t+yztx^2+tyzx^2,\text{ i.e. }d_3(y,z,t;1,x^2)=S_3(y,z,t)x^2.
\]
This together with (\ref{(16)}), (\ref{(17)}) and (\ref{(18)}) implies that 
\[
d_3(y,z,t;x^2,1)=0=x^2S_3(y,z,t)
\]
is a consequence of $xy^2x=0$.
It follows from $xyx^2=0$ that $xyzx^2=yx^2zx=0$ and we can apply similar arguments.

4) Obviously $x^2y^2x=0$ follows from $xy^2x=0$.

We linearize in $y$ and replace $yx$ instead of $y$. Then we replace $zx$ instead of $z$ and, subtracting these two identities we obtain $x^2[y,z]x=0$.

It is obvious that $x^2[y,z]x=0$ and $x^2y^2x=0$ follow from $x^2yx=0$.

5) Now we shall obtain that
\[
2d_3(y,z,t;x^2,1)+x^2S_3(y,z,t)=0
\]
follows from
\begin{equation}\label{(19)}
x^2[y,z]-[y,z]x^2=0.
\end{equation}

Using (\ref{(13)}) and (\ref{(14)}) we obtain from (\ref{(19)})
\[
tx^2yz-2zx^2yt+2yx^2tz-yx^2zt=0.
\]
In (\ref{(19)}) we replace $z$ by $z\circ t=zt+tz$ and, taking into account (\ref{(12)}) we find that
\[
x^2(z\circ t)y+y(z\circ t)x^2=0.
\]
This and (\ref{(14)}) imply
\[
yx^2tz-zx^2yt=-yx^2zt+tx^2yz.
\]
Together with (\ref{(19)}) the latter equality gives
\begin{equation}\label{(20)}
tx^2yz=yx^2zt
\end{equation}
and, similarly, $yx^2zt=zx^2ty$. We multiply (\ref{(19)}) from the left by $t$ and obtain from (\ref{(20)})
\begin{equation}\label{(21)}
2x^2tyz=tx^2[z,y].
\end{equation}
In virtue of (\ref{(14)}) and (\ref{(20)}) we obtain $x^2yzt+yztx^2=0$. Hence, it can be seen from (\ref{(12)})
that all monomials which start or end with $x^2$ are equal up to the sign. Hence if we multiply (\ref{(19)}) by $x$ from the left
we obtain $tx^2[y,z]+4tyzx^2=0$. From here and from (\ref{(15)}) we derive that $tx^2[y,z]=0$ and $tyzx^2=0$.
Hence $d_3(y,z,t;x^2,1)=0$ and $S_3(y,z,t)x^2=0$. In particular, this gives that $x^2[y,z]x=0$ and $x^2y^2x=0$ are also
consequences of (\ref{(19)}).

6) Now we shall obtain the consequences of
\begin{equation}\label{(22)}
yx^2z-zx^2y=0.
\end{equation}
Replacing $z$ by $zx$ and $y$ by $yx$, respectively, and subtracting these equalities, we obtain $x^2[y,z]x=0$.
If we replace $z$ by $yx$ in (\ref{(22)}) we shall obtain $x^2y^2x=0$. Replacing $z$ by $zt$ in (\ref{(22)})
and taking the alternating sum in $y,z,t$ we shall obtain
\begin{equation}\label{(23)}
d_3(y,z,t;x^2,1)-d_3(y,z,t;1,x^2)=0.
\end{equation}
Together with (\ref{(16)}) this implies $d_3(y,z,t;x^2,1)=0$. Now we linearize (\ref{(22)}), multiply it from the left by $x$ and
taking into account (\ref{(15)}) we derive
\[
yx^2tz+2ytzx^2+ytx^2z-ztyx^2-zx^2ty=0.
\]
The alternating sum in $y,z,t$ gives
\[
2d_3(y,z,t;x^2,1)+4S_3(y,z,t)x^2+2d_3(y,z,t;1,x^2)=0.
\]
Then (\ref{(16)}) implies that $S_3(y,z,t)x^2=0$.

7) Now we shall consider the consequences of
\begin{equation}\label{(24)}
\alpha(x^2[y,z]-[y,z]x^2)+\beta(yx^2z-zx^2y)=0,\,\alpha\not=0,\beta\not=0.
\end{equation}
We multiply from the right by $t$ and take the alternating sum in $y,z,t$:
\begin{equation}\label{(25)}
2x^2S_3(y,z,t)+(\alpha+\beta)d_3(y,z,t;x^2,1)=0.
\end{equation}
Similarly, multiplying from the left by $t$ and taking into account (\ref{(25)}), we obtain
\[
d_3(y,z,t;x^2,1)=x^2S_3(y,z,t)=0.
\]

Now we multiply (\ref{(24)}) from the right by $x$:
$(\alpha-\beta)x^2[y,z]x=0$.
If $\alpha=\beta$ we can check explicitly that again $x^2[y,z]x=0$.

If we replace $z$ by $yx$ in (\ref{(24)}) we shall obtain $(3\alpha-\beta)x^2y^2x=0$.
If $3\alpha=\beta$ we replace $z$ by $xy$ in (\ref{(24)}).

The case in Fig. 6 (c) of the theorem is a consequence of the considered cases 3)--7).

In this way, the proof of the theorem is completely finished.

\end{document}